\documentclass[11pt]{article}
\usepackage{amsmath,amsfonts,amssymb,amsthm}
 \def\N{\mathbb{ N}}
\newtheorem{thm}{Theorem}
\newtheorem{cor}{Corollary}
\newtheorem{prop}{Proposition}
\newtheorem{definition}{Definition}

\newenvironment{remark}{\smallskip\noindent{\bf Remark. }}{\smallskip}
\newcommand\bi[2]{{{#1}\atopwithdelims(){#2}}}
\newcommand\qbi[3]{{{#1}\atopwithdelims[]{#2}}_{#3}}
\begin{document}

\begin{center}
{\Large\Large\bf Two new families of $q$-positive
integers}
\end{center}

\vskip 2mm \centerline{Sharon J. X.
Hou$^{a,}$\,\footnote{corresponding author.} and Jiang
Zeng$^{a,\,b}$}

\begin{center} $^a$Center for Combinatorics, LPMC\\[5pt]
Nankai University, Tianjin 300071, People's Republic of China\\[5pt]
{\tt jxhou@eyou.com}
\end{center}

\begin{center}
$^b$Institut Girard Desargues,
Universit\'e Claude Bernard (Lyon I)\\[5pt]
F-69622 Villeurbanne Cedex, France \\[5pt]
{\tt zeng@desargues.univ-lyon1.fr}
\end{center}
{\small \vskip 0.7cm \noindent{\bf Abstract.} Let $n,p,k$ be three
positive integers. We prove that the rational fractions of $q$:
$${n \brack k}_{q}
{}_3\phi_{2}
\left[\left.\begin{matrix}q^{1-k},q^{-p},q^{p-n}\\[5pt]
q,q^{1-n}\end{matrix}\right|
q;q^{k+1}\right]\quad\textrm{and}\quad q^{(n-p)p}\qbi{n}{k}{q} \
{}_3\phi_2\!\left[\left.\begin{matrix}q^{1-k},q^{-p},q^{p-n}\\[5pt]
q,q^{1-n}\end{matrix}\right|q;q\right]
$$
are polynomials of $q$ with positive integer coefficients.  This
generalizes a recent result of Lassalle (Ann. Comb. 6(2002), no.
3-4, 399-405), in the same way as the classical $q$-binomial
coefficients refine the ordinary binomial coefficients. }


\vskip 0.2mm
\noindent{\bf MR Subject Classifications}: Primary 05A30; Secondary 33D15;

\section{Introduction}
In~\cite{La} Lassalle introduced a new family of positive integers
generalizing the classical binomial coefficients. We refer the
readers to \cite{La2} for  some motivations from the study of Jack
polynomials and to \cite{JLZ} for some further results and
extensions of these coefficients. In this paper we will present a
different generalization of Lassalle's coefficients. More
precisely, using basic hypergeometric functions we will show that
all the results of Lassalle~\cite{La} have {\em natural}
$q$-analogues.

In order to present its $q$-analogues we need first to introduce
some notations. Throughout this paper $q$ is a complex variable
such that $|q|<1$. We use the standard notation of \cite{GR} for
basic hypergeometric series
\[{}_{r}\phi_s\!\left[\left.\begin{matrix}a_1,a_2,\dots,a_r\\[5pt]
b_1,b_2,\dots,b_s\end{matrix}\right| q;z\right]= \sum _{k=0}
^{\infty}\frac {(a_1,\dots, a_{r}; q)_k} {(b_1,\dots, b_s;
q)_k}(-1)^{(1+s-r)k}q^{(1+s-r){k\choose 2}}\frac{z^k}{(q;q)_k},\]
where for an indeterminate $a$ and some positive integer $k$, the
$q$-raising factorial is defined by ${(a;q)}_k = (1-a)(1-aq)
\ldots (1-aq^{k-1})$ and $(b_1,\dots, b_s;
q)_k=(b_1;q)_k\ldots(b_s;q)_k$. The $q$-binomial coefficients
${n\brack k}_q$ are usually defined by
 $${n\brack k}_q =
\frac{(q;q)_n}{(q;q)_k (q;q)_{n-k}},
$$
which have two less obvious definitions as follows:
\begin{equation}\label{eq:qbin}
\qbi{n}{p}{q} = {}_2\phi_1\!\left[\left.\begin{matrix}
q^{-p},q^{p-n}\\[5pt]
q\end{matrix}\right|q;q^{n+1}\right]=q^{p(p-n)}
{}_2\phi_1\!\left[\left.\begin{matrix}
q^{-p},q^{p-n}\\[5pt]
q\end{matrix}\right|q;q\right].
\end{equation}
Actually these expressions are obtained by specializing
$a=q^{p-n}, c=q$ in the celebrated $q$-Chu-Vandermonde
formula~\cite[p. 236]{GR}:
\[\frac{(c/a;q)_p}{(c;q)_p} = {}_2\phi_1\!\left[\left.
\begin{matrix}q^{-p},a\\[5pt]
c\end{matrix}\right|q;cq^p/a \right]=a^{-p} {}_2\phi_1\!
\left[\left.
\begin{matrix}q^{-p},a\\[5pt]
c\end{matrix}\right|q;q\right].
\]
Of course using this  relation as a definition of $q$-binomial
coefficients would be rather tautological. However, quite
surprisingly, it is possible to define two new families of {\em
$q$-positive integers}, i.e., a polynomial of $q$ with {\em non
negative integer coefficients}, by slightly modifying the
$q$-Chu-Vandermonde formula. In fact there are two such
$q$-analogues of Lassalle's generalized binomial coefficients. In
the next two sections we present the  first $q$-analogue and its
{\em raison d'\^etre} in  the context of linearization problem of
$q$-binomial coefficients. In Section~4 we outline the second
$q$-analogue.
 We end this paper with the first values
of our generalized $q$-binomial coefficients., which seems to
suggest some \emph{unimodal} properties of the coefficients of
these polynomials.

\section{The first $q$-analogue}
\begin{definition}
For any positive integers $n,p,k$, define
\begin{eqnarray}
\bi{n}{p,k}_{q}&=& q^{(n-p)p}\qbi{n}{k}{q} \
{}_3\phi_2\!\left[\left.\begin{matrix}q^{1-k},q^{-p},q^{p-n}\\[5pt]
q,q^{1-n}\end{matrix}\right|q;q\right].\label{eq:def1}
\end{eqnarray}
\end{definition}
We have  obviously
\[\bi{n}{p,k}_{q} = 0 \quad \textrm{for} \quad
k>n \quad \textrm{or} \quad p>n,
\]
and for $p,k\leq n$,
\[
 \qquad \bi{n}{p,k}_{q}
=\bi{n}{n-p,k}_{q},\qquad  \bi{n}{p,n}_{q}=\qbi{n}{p}{q},
\]
the last equations following directly from the $q$-Chu-Vandermonde
formula.

Set $[n]_q=1+q+\cdots +q^{n-1}=(1-q^n)/(1-q)$ for $n\geq 0$. We
can rewrite \eqref{eq:def1} as follows:
\begin{eqnarray}\label{qdef}
\bi{n}{p,k}_{q} &=&\frac{[n]_q}{[k]_q}\sum_{r \ge 0}
\qbi{p}{r}{q}\qbi{n-p}{r}{q}
\qbi{n-r-1}{k-r-1}{q}q^{(n-p)p+r(r-k)}.\label{eq4}
\end{eqnarray}
Therefore
\begin{equation*}
\bi{n}{0,k}_{q} = \qbi{n}{k}{q}, \qquad \bi{n}{1, k}_{q} = [k]_q
q^{n-k} \qbi{n}{k}{q},
\end{equation*}
and
\begin{eqnarray*}
\bi{n}{2,k}_{q} &=& q^{2n-3-k}[k]_q \qbi{n}{k}{q}
+q^{2(n-k)}\frac{[n]_q[n-3]_q}{[2]_q} \qbi{n-2}{k-2}{q}.
\end{eqnarray*}
For $p>0$, we have also
\[
\bi{n}{p, 0}_{q} = 0, \quad \bi{n}{p, 1}_{q} = q^{(n-p)p}[n]_q,\]
\[\bi{n}{p, 2}_{q}=q^{(n-p)p}\frac{[n]_q}{[2]_q}
\left([n-1]_q+[p]_q[n-p]_q\right).
\]
These results suggest that $\bi{n}{p, k}{q}$ are polynomials of
$q$ with non negative integer coefficients. Indeed, it's not hard
to see that they are polynomials with integer coefficients. Recall
the following  two transformations of Sears~\cite[p. 61]{GR}:
\begin{eqnarray}\label{q-whip-2}
_3\phi_2\!\left[\left.\begin{matrix}q^{-n},a,b\\[5pt]c,d\end{matrix}\right|q;cdq^n/ab\right]=
\frac{(c/a;q)_n}{(c;q)_n}\!\
_3\phi_2\left[\left.\begin{matrix}q^{-n},a,d/b\\[5pt]d,q^{1-n}a/c\end{matrix}\right|q;q
\right],
\end{eqnarray}
and
\begin{eqnarray}\label{q-whip-1}
_3\phi_2\!\left[\left.\begin{matrix}q^{-n},a,b\\[5pt]c,d\end{matrix}\right|q;q\right]=
\frac{(c/a;q)_n}{(c;q)_n}a^n\!\
_3\phi_2\left[\left.\begin{matrix}q^{-n},a,d/b\\[5pt]d,q^{1-n}a/c\end{matrix}\right|q;bq/c
\right].
\end{eqnarray}
Applying \eqref{q-whip-1} to \eqref{eq:def1}, we get
\begin {eqnarray}
\bi{n}{p,k}_{q} &=& \qbi{n}{k}{q}\qbi{n}{p}{q}
{}_3\phi_2\!\left[\left.\begin{matrix}q^{-p},q^{p-n}, q^{k-n}\\[5pt]
q^{1-n}, q^{-n}\end{matrix}\right|q;q^{1-k}\right]\nonumber\\[5pt]
&=&\sum_{i \geq 0} (-1)^i
\qbi{n-i}{k}{q}\frac{[n]}{[n-i]}\qbi{n-i}{i}{q} \qbi{n-2i}{p-i}{q}
q^{i(i-1)/2}.\label{eq:expression2}
\end {eqnarray}
Since
\begin{equation*}
\frac{[n]}{[n-i]}\qbi{n-i}{i}{q}=\frac{1-q^{n-i}+q^{n-i}-q^n}{1-q^{n-i}}
\qbi{n-i}{i}{q}=\qbi{n-i}{i}{q}+q^{n-i}\qbi{n-i-1}{i-1}{q},
\end{equation*}
we derive that $\bi{n}{p,k}_{q}$ is $q$-integer, i.e.,
polynomials of $q$ with integer coefficients.

Two generating functions for ${n \choose p,k}_1$  were published
in \cite{La}. We consider  their $q$-analogues below.
\begin{thm} There holds
\begin{eqnarray}
& &\sum_{0 \leq k,p \leq n}{n \choose  p,k}_qx^p q^{p \choose 2}
y^k q^{k \choose 2} \nonumber\\[5pt]
&=& (-y;q)_n(-x;q)_n\sum_{i \geq 0} \frac{1-q^n}{1-q^{n-i}}{n-i
\brack
i}_q\frac{(-x)^iq^{i^2-i}}{(-yq^{n-i};q)_i(-x;q)_{i}(-xq^{n-i};q)_i}.\label{gf1}
\end{eqnarray}
\end{thm}
\begin{proof} Using \eqref{eq:expression2} one has
\begin{eqnarray*}
& &\sum_{0 \leq k,p \leq n}{n \choose p,k}_qx^p q^{p \choose 2} y^k q^{k \choose 2} \\[5pt]
& = & \sum_{0 \leq k,p, \leq n}\sum_{i \geq 0}^n(-1)^i
\frac{1-q^n} {1-q^{n-i}}{n-i \brack k}_q{n-i \brack i}_q{n-2i
\brack p-i}_qq^{{i \choose 2}}y^k q^{k \choose 2} x^p q^{p\choose 2}\\[5pt]
& = & \sum_{i \geq 0}(-1)^i \frac{1-q^n}{1-q^{n-i}}x^i \sum_{k
\geq 0}{n-i \brack k}_q y^k q^{k \choose 2} \sum_{p \geq
0}{n-2i \brack  p-i}_q(xq^{i})^{p-i}q^{p-i \choose 2}q^{i^2-i}\\[5pt]
& = & \sum_{i \geq 0}(-1)^i \frac{1-q^n}{1-q^{n-i}}x^i
(-y;q)_{n-i} (-xq^{i};q)_{n-2i}q^{i^2-i},
\end{eqnarray*}
which yields \eqref{gf1}.
\end{proof}

\begin{remark} When $q=1$, since
$$
\sum_{k<n}\bi{n-k}{k}\frac{n}{n-k}z^k=\left(\frac{1+\sqrt{1+4z}}{2}\right)^n
+\left(\frac{1-\sqrt{1+4z}}{2}\right)^n,
$$
setting $z=\frac{-x}{(1+x)^2(1+y)}$, the above generating
functions  can be written as follows:
\[
\sum_{0 \leq k,p \leq n}{n \brack p,k}_1x^p
y^k=2^{-n}[(1+y)(1+x)]^n \left((1+\sqrt{1+4z})^n
+(1-\sqrt{1+4z})^n\right).
\]
\end{remark}

The following is the $q$-analogue of Lassalle's recurrence
relation in \cite[Lemma~3.2]{La}.
\begin{prop}
For $k \neq 0$ and $0 \leq p \leq n$, we have
\begin{eqnarray}
(1-q^{n-p+1}){n \choose p-1, k}_q -(1-q^p){n \choose p, k}_q =
\frac{[n]_q}{[n-1]_q}q^{p-1}(1-q^{n-2p+1}){n-1 \choose p-1,
k}_q.\label{conti2}
\end{eqnarray}
\end{prop}
\begin{proof}
Indeed, up to the factor $[n]_q/[k]_q$ and using \eqref{qdef}, the
left-hand side can be written as
\begin{eqnarray*}
& &\sum_{r \geq 0}q^{A}{n-r-1 \brack k-r-1}_q
\left((1-q^{n-p+1})q^{p-1}{p-1 \brack
r }_q{n-p+1 \brack r}_q-(1-q^p)q^{n-p}{p \brack r}_q{n-p \brack r}_q\right) \\[5pt]
& = &(1-q^{n-2p+1})\sum_{r \geq 0}q^{A}{n-r-1 \brack k-r-1}_q
\left(q^{p-1}{p-1 \brack r}_q{n-p \brack r}_q
- q^{n-2r+1}{p-1 \brack r-1}_q{n-p \brack r-1}_q\right)\\[5pt]
&=& (1-q^{n-2p+1})\sum_{r \geq 0}q^{A}{p-1 \brack r}_q{n-p \brack
r}_q q^{p-1}\left({n-r-1 \brack k-r-1}_q -q^{n-k}
{n-r-2 \brack k-r-2}_q\right)\\[5pt]
&=& (1-q^{n-2p+1})\sum_{r \geq 0}q^{A}{p-1 \brack r}_q{n-p \brack
r}_q{n-r-2 \brack k-r-1}_qq^{p-1},
\end{eqnarray*}
where $A=r(r-k)+(n-p)(p-1)$. The result follows then from Equation
\eqref{eq4}.
\end{proof}

For the second generating function, we will need the $q$-binomial
formula~\cite[236]{GR}:
\begin{equation}\label{eq:q-bi}
\sum_{n\geq
0}\frac{(a;q)_n}{(q;q)_n}x^n=\frac{(ax;q)_\infty}{(x;q)_\infty},
\end{equation}
and Jackson's $q$-Pfaff-Kummer transformation~\cite[p.241]{GR}:
\begin{eqnarray}\label{eq:jac}
 {}_2\phi_1\!
\left[\left.\begin{matrix} a,b\\[5pt]
c\end{matrix}\right|q;z\right] =\frac{(az;q)_\infty}{(z;q)_\infty}
{}_{2}\phi_2\!\left[\left.\begin{matrix}a,c/b\\[5pt]
c, az\end{matrix}\right| q;bz\right].
\end{eqnarray}

\begin{thm}
For $0 \leq p \leq n$, we have
\begin{eqnarray}
\sum_{k \geq 1}{n \choose p,k}_qy^k q^{{k \choose
2}}&=&yq^{p(n-p)}(-y;q)_p[n]_q\,
{}_{2}\phi_1\!\left[\left.\begin{matrix}q^{p+1},q^{p-n+1}\\[5pt]
q^{2} \end{matrix}\right| q; -yq^{n-p}\right].\label{eq:gf2}
\end{eqnarray}
\end{thm}
\begin{proof} By applying \eqref{q-whip-1} to \eqref{eq:def1} and then applying
\eqref{q-whip-2} we get:
\begin{eqnarray}
{n \choose p,k}_q&=& q^{(n-p)p}{n \brack k}\!\
\frac{(q^{1-p};q)_{k-1}}{(q^{1-n};q)_{k-1}}q^{(k-1)(p-n)}
{}_3\phi_2\!\left[\left.\begin{matrix}q^{1-k},q^{p-n},q^{1+p}
\\[5pt]q^{p-k+1};q\end{matrix}\right|q;q^{n-p} \right]\nonumber\\[5pt]
&=&q^{(n-p)p}{n \brack
k}\frac{(q^{1-p};q)_{k-1}}{(q^{1-n};q)_{k-1}}q^{(p-n)(k-1)}
\frac{(q^{-k};q)_{k-1}}
{(q^{p-k+1};q)_{k-1}}{}_3\phi_2\!\left[\left.\begin{matrix}q^{1-k},q^{1+p},q^{n-p+1}
\nonumber\\[5pt]
q, q^2
\end{matrix}\right|q;q\right] \\[5pt]
&=& (-1)^{k-1}q^{(n-p)p-\frac{(k-1)(k+2)}{2}}[n]_q
{}_3\phi_2\!\left[\left.\begin{matrix}q^{1-k},q^{1+p},q^{n-p+1}\label{def2bis}\\[5pt]
q, q^2
\end{matrix}\right|q;q\right].
\end{eqnarray}
Hence
\begin{eqnarray*}
\sum_{k \geq 1}{n \choose p,k}_qy^kq^{k \choose 2 }
&=&yq^{p(n-p)}[n]_q\sum_{l \geq 0}\frac{(q^{p+1};q)_l
(q^{n-p+1};q)_l }{(q^2;q)_l(q;q)_l(q;q)_l} q^l \sum_{k \geq
l}(q^{-k};q)_l(-yq^{-1})^k.
\end{eqnarray*}
But the $q$-binomial formula \eqref{eq:q-bi} implies that
\[
\sum_{k \geq l}(q^{-k};q)_l(-yq^{-1})^k = (-1)^lq^{{l\choose
2}}\sum_{k\geq l}(q^{k-l+1};q)_l(-yq^{-l-1})^k = \frac{(q;q)_l}{(1+yq^{-1})(-y^{-1}q^2;q)_l}.
\]
Hence
\begin{equation}
\sum_{k \geq 1}{n \choose  p,k}_qy^kq^{{k \choose
2}}=\frac{yq^{p(n-p)}}{1+yq^{-1}}[n]_q \sum_{l\geq
0}\frac{(q^{p+1};q)_l(q^{n-p+1};q)_lq^l}{(-y^{-1}q^2;q)_l(q^2;q)_l(q;q)_l}.
\end{equation}
Using the formula
\begin{eqnarray}\label{q-inverse}
(a;q)_n=(a^{-1};q^{-1})_n(-a)^nq^{n\choose 2} \end{eqnarray}
 and Jackson's transformation, we can rewrite the above sum as follows:
\begin{eqnarray*}
\sum_{k \geq 1}{n \choose  p,k}_qy^kq^{{k \choose 2}}
&=&\frac{yq^{p(n-p)}}{1+yq^{-1}}[n]_q
{}_{2}\phi_2\!\left[\left.\begin{matrix}q^{-1-p},q^{p-n-1}\\[5pt]
q^{-2}, -yq^{-2}\end{matrix}\right| q^{-1}; -yq^{n-2}\right]\\[5pt]
&=& yq^{p(n-p)}(-y;q)_p[n]_q
{}_{2}\phi_1\!\left[\left.\begin{matrix}q^{-1-p},q^{p-n-1}\\[5pt]
q^{-2} \end{matrix}\right| q^{-1}; -yq^{p-1}\right].
\end{eqnarray*}
We then recover  \eqref{eq:gf2} by applying \eqref{q-inverse}
again.
\end{proof}

\noindent{\bf Remark}. When $q=1$, the above theorem reduces to
Lassalle's generating function~\cite{La}, which was proved by
using induction and contiguous relation. One could give another
proof of Theorem 2  by using Proposition 1.

The following result is crucial to prove that ${n \choose p, k}_q$
is a $q$-positive integer.

\begin{cor} We have
\[ {n \choose p,k}_q=\frac{[n]_q}{[p]_q}\,\sum_{i=0}^{k-1}{n-p
\brack k-1-i}_q{n-p+i \brack i}_q{p \brack i+1}_q
q^{(i+1)(i+1-k)+(n-p)p}.
\]
\end{cor}
\begin{proof}
Using the $q$-binomial formula we have
\[(-y; q)_p = \sum_{j=0}^p{p \brack j}y^jq^{{j \choose 2}}.
\]
Extracting the
coefficient of $y^k$ in \eqref{eq:gf2} we obtain
$$
{n \choose p,k}_q=\frac{[n]_q}{[n-p]_q}\sum_{i=0}^{k-1} {p \brack
k-1-i}_q{p+i \brack i}_q{n-p \brack i+1}_qq^{(n-p)p+(i+1)(i+1-k)}.
$$
As ${n \choose p,k}_q={n \choose n-p,k}_q$, substituting $p$ by
$n-p$ yields the desired identity.
\end{proof}
By applying
\begin{eqnarray*}
{n \choose p,k}_q&=& \frac{[p]_q}{[n]_q}{n \choose
p,k}_q+q^p\frac{[n-p]_q}{[n]_q}{n \choose n-p,k}_q,
\label{q-positive}
\end{eqnarray*}
we can write
\begin{eqnarray*}
{n \choose p,k}_q &=& \sum_{i=0}^{k-1}{n-p \brack
k-1-i}_q{n-p+i \brack i}_q{p \brack i+1}_q q^{(n-p)p+(i+1)(i+1-k)} \\[5pt]
&+& \sum_{i=0}^{k-1}{p \brack k-1-i}_q {p+i \brack i}_q
{n-p \brack i+1}_q q^{(n-p)(p+1)+(i+1)(i+1-k)},
\end{eqnarray*}
which implies the following
\begin {thm}
The polynomials ${n \choose p,k}_q$ are $q$-positive integers.
\end{thm}

Since $q$-binomial coefficients have various nice combinatorial
interpretations, it would be possible to derive a combinatorial
interpretation for ${n \choose p,k}_q$ from the above expression.

\section{Further extensions}
It is surprising that the general numbers~$c_k^{({\bf r})}$ of
Lassalle~\cite{La} have also a $q$-analogue, which are also
$q$-positive integers. We shall explain such a $q$-analog in this
section. Note that the $q$-Chu-Vandermonde formula:
\begin{equation}\label{eq:va}
{n+m\brack k}_q = \sum_{i\ge 0} {n\brack i}_q {m\brack k-i}_q
q^{(n-i)(k-i)}
\end{equation}
implies that
\begin{eqnarray}
{x\brack r_1}_q {x\brack r_2}_q &=& {x\brack r_1}_q \sum_{k \ge 0}
{r_1\brack k}_q {x-r_1\brack r_2-k}_q q^{(r_1-k)(r_2-k)}
\nonumber\\[5pt]
&=& \sum_{k \ge 0} q^{(r_1-k)(r_2-k)} {r_1+r_2-k\brack
k,r_1-k,r_2-k}_q {x\brack r_1+r_2-k}_q. \label{eq:qbinom1}
\end{eqnarray}
Set ${\bf r}=(r_1,\ldots, r_m)$ and $|{\bf r}|=r_1+\cdots +r_m$.
Iterating \eqref{eq:qbinom1} yields:
\begin{equation}\label{eq:qlinbin}
{x\brack r_1}_q \cdots {x\brack r_m}_q = \sum_{l\geq 0} d_l({\bf
r};q)  {x\brack l}_q,
\end{equation}
where $d_l({\bf r};q)$ are $q$-positive integers given by

\begin{eqnarray}
d_l({\bf r};q)&=&\sum_{k_1,\ldots, k_{m-2}\geq 0}
{r_1+r_2-k_1\brack k_1,r_1-k_1, r_2-k_1}_q
{r_1+r_2+r_3-k_1-k_2\brack k_2,r_1+r_2-k_1-k_2, r_3-k_2}_q \nonumber\\[5pt]
 & &\times \cdots \times
 {r_1+\ldots +r_{m-2}-k_1-\ldots -k_{m-3}\brack k_{m-3},r_1+\ldots +r_{m-3}-k_1-\ldots
 -k_{m-3}, r_{m-2}-k_{m-3}}_q \nonumber\\[5pt]
& &\times {l\brack r_m}_q{r_m\brack |{\bf r}|-k_1-\ldots-k_{m-2}-l
}_q q^B,\label{dl-rm}
\end{eqnarray}
where
\begin{eqnarray*}
B&=&(r_1-k_1)(r_2-k_1)+(r_1+r_2-k_1-k_2)(r_3-k_2)+\cdots \\[5pt]
&&+ (r_1+r_2+\cdots +r_{m-2}-k_1-\cdots
-k_{m-2})(r_{m-1}-k_{m-2})\\[5pt]
&&+(l-r_m)(l-r_1-\cdots -r_{m-1}+k_1+\cdots+ k_{m-2}).
\end{eqnarray*}
In particular, for $m=2$ we have
\begin{equation}\label{m=2}
d_l(r_1,r_2;q)=q^{(l-r_1)(l-r_2)} {l\brack r_1}_q {r_1\brack
l-r_2}_q.
\end{equation}
On the other hand, identity (\ref{eq:va}) implies also
\begin{equation}\label{eq:qbinom2}
 {x+r_1-1\brack r_1}_q = \sum_{k
\ge 0} {r_1-1\brack r_1-k}_q {x\brack k}_q q^{k(k-1)}.
\end{equation}
From \eqref{eq:qlinbin} and \eqref{eq:qbinom2} we derive the
following result.

\begin{equation*}\label{eq:qlinlas}
{x+r_1-1\brack r_1}_q\cdots {x+r_m-1\brack r_m}_q = \sum_{l\geq 0}
\tilde{c}_l({\bf r};q) {x\brack l}_q.
\end{equation*}
where
\begin{equation}
\tilde{c}_l({\bf r};q)=\sum_{\bf k}d_l({\bf
k};q)\prod_{i=1}^m{r_i-1\brack k_i-1}_q q^{k_i(k_i-1)}.
\end{equation}

\begin{thm}
The polynomial
$$
c_l({\bf r};q)=\frac{[r_1+\ldots + r_m]_q}{[l]_q}\tilde{c}_l({\bf
r};q)
$$
is a $q$-positive integer.
\end{thm}
\begin{proof} By \eqref{dl-rm} there is a polynomial
$P_m({\bf k};q)\in \N[q]$ such that
$$
d_l({\bf k};q)=\frac{[l]_q}{[k_m]_q}P_m(k; q).
$$
Since $d_l({\bf k};q)$ is symmetric with respect to ${\bf
k}=(k_1,\ldots, k_m)$, the above formula infers that there is a
polynomial $P_i({\bf k};q)\in \N[q]$ for each $j\in \{1,\ldots,
m\}$ such that
$$
d_l({\bf k};q)=\frac{[l]_q}{[k_j]_q}P_j({\bf k};q).
$$
Therefore, using \eqref{eq:qlinlas} we have
\begin{eqnarray*}
c_l({\bf r};q)&=&\sum_{j=1}^m\frac{[r_j]q^{r_1+\cdots
+r_{j-1}}}{[r_1+\cdots +r_m]}\, c_l({\bf r};q)\\[5pt]
&=&\sum_{j=1}^mq^{r_1+\cdots +r_{j-1}} \sum_{\bf k}P_i({\bf
k};q){r_j\brack k_j}_qq^{k_j(k_j-1)}\prod_{i\neq j}{r_i-1\brack
k_i-1}_qq^{k_i(k_i-1)},
\end{eqnarray*}
which is clearly a $q$-positive integer.
\end{proof}

We can also derive a  simpler  formula for
$\tilde{c}_j(r_1,\dots,r_m;q)$ using the {\em $q$-difference
operator}. Set $[x;q]={(q^x-1)/(q-1)}$ and
$$
[x;q]_n=[x;q][x-1;q]\cdots [x-n+1;q]=\frac{(q^{x-n+1};q)_n}{(1-q)^n}.
$$
 We define  the
$q$-difference operator $\Delta_q$ by
$$
\Delta_q^0f(x)=f(x),\, \Delta_q^{n+1}f(x)=\Delta_q^n(E-q^nI)f(x),
$$
where $If(x)=f(x)$ and $Ef(x)=f(x+1)$. Note that
$$
\Delta_q^nf(x)=(E-q^{n-1}I)(E-q^{n-2}I)\cdots (E-I)f(x).
$$
By the  $q$-Chu-Vandermonde formula, we have
\begin{equation}\label{eq:qdiff}
\Delta_q^nf(x)=\sum_{k=0}^n(-1)^k{n\brack k}_q q^{k\choose
2}f(x+n-k).
\end{equation}

It's easy to see that
$$
\Delta_q^n[x;q]_m=[m;q]_n[x;q]_{m-n}q^{n(x+n-m)}.
$$
We have also
\begin{equation}\label{eq:qnewton}
p(x)=\sum_{n\geq 0}\frac{\Delta_q^np(0)}{[n]!}[x;q]_n.
\end{equation}
It follows from \eqref{eq:qnewton}, \eqref{eq:qdiff} and
\eqref{eq:qlinlas} that
\begin{equation}\label{qdicoeff}
\tilde{c}_k(r_1,\dots,r_m;q)=\sum_{j=1}^k(-1)^{k-j}\qbi{k}{j}{q}q^{{k-j\choose
2}}\prod_{l=1}^m\qbi{j+r_l-1}{r_l}{q}.
\end{equation}
Set ${\bf r}=(r_1,\dots,r_m)$ and
\begin{eqnarray}\label{eq:q1}
c_k({\bf r};q)&=&\frac{[r_1+\ldots +
r_m]_q}{[k]_q}\tilde{c}_k(r_1,\dots,r_m;q)\nonumber\\[5pt]
&=&\sum_{i=1}^mq^{r_1+\cdots +r_{i-1}}\frac{[r_i]_q}{[k]_q}
\sum_{j=1}^k(-1)^{k-j}\qbi{k}{j}{q}q^{{k-j\choose
2}}\prod_{l=1}^m\qbi{j+r_l-1}{r_l}{q}\nonumber\\[5pt]
&=& \sum_{i=1}^m\sum_{j=1}^k(-1)^{k-j}q^{r_1+\cdots
+r_{i-1}+{k-j\choose 2}}\nonumber\\[5pt]
&&\hspace{2cm}\times\qbi{k-1}{j-1}{q}\qbi{j+r_i-1}{r_i}{q}\prod_{l=1,\,l\neq
j}^m\qbi{r_l+i-1}{r_l}{q}.
\end{eqnarray}
Thus we have obtained another proof of the $q$-integrality of
$c_k({\bf r};q)$.

Finally, when $\mathcal{\mathbf{r}}=(r_1,r_2)$, we have the
following result.
\begin{thm} The coefficients  ${r_1+r_2 \choose
r_1,k}_q$ satisfy
\begin{eqnarray}\label{eq:qlinlas2}
{x+r_1-1 \brack r_1}_q{x+r_2-1 \brack r_2 }_q= \sum_{k \geq
1}\frac{[k]_qq^{k(k-1)-r_1r_2}}{[r_1+r_2]_q}{r_1+r_2 \choose
r_1,k}_q{x \brack k}_q.
\end{eqnarray}
\end{thm}

\begin{proof}
By  \eqref{qdicoeff} we have
\begin{eqnarray*}
c_k(r_1,r_2;q) &=&\frac{[r_1+r_2]_q}{[k]_q} \sum_{j=
0}^k(-1)^{k-j}\qbi{k}{j}{q}q^{{k-j\choose 2}}
\qbi{j+r_1-1}{r_1}{q}\qbi{j+r_2-1}{r_2}{q}\\[5pt]
&=& (-1)^{k-1}[r_1+r_2]_qq^{{k-1\choose 2}}
{}_3\phi_2\!\left[\left.\begin{matrix}q^{1-k},q^{1+r_1},q^{1+r_2}\\[5pt]
q, q^2
\end{matrix}\right|q;q\right].
\end{eqnarray*}
Comparing with \eqref{def2bis} we see that
$c_k(r_1,r_2;q)=q^{k(k-1)-r_1r_2}{r_1+r_2 \choose r_1,k}_q$.
\end{proof}

\section{The Second $q$-analogue}
In this section, we give another new family of q-positive integers
${n \brack p,k}_q$, which have the similar properties as ${n
\choose p,k}_q$. Since the proofs are similar we omit the details.
\begin{definition}
For any positive integers $n,p,k$, define \begin{eqnarray}
 {n \brack p, k}_q = {n \brack k}_{q} \, _3\phi_2\!\left[\left.\begin{matrix}q^{1-k},q^{-p},q^{p-n}\\[5pt]
q,q^{1-n}\end{matrix}\right|q;q^{k+1}\right].\label{eq:def2}
\end{eqnarray}
\end{definition}
Obviously we have
$$\qbi{n}{p,k}{q} =\qbi{n}{n-p,k}{q},
\qbi{n}{p,n}{q}=\qbi{n}{p}{q}.$$ The definition \eqref{eq:def2}
could be written as follows:
$$\qbi{n}{p,k}{q}=\frac{[n]_q}{[k]_q}\sum_{r \ge 0}
\qbi{p}{r}{q}\qbi{n-p}{r}{q}
\qbi{n-r-1}{k-r-1}{q}q^{r^2}.\label{eq3}
$$

Therefore
\[{n \brack 0,k}_q ={n \brack k}_q, \quad {n, \brack 1,k}_q
=[k]_q{n \brack k}_q,\]
 and
\[\qbi{n}{2,k}{q} = [k]_q \qbi{n}{k}{q} +q^{2}
\frac{[n]_q[n-3]_q}{[2]_q} \qbi{n-2}{k-2}{q}.\\[5pt]
\]

For $p \geq 0$, we have also
\[ \qbi{n}{p, 0}{q} =0, \quad \qbi{n}{p, 1}{q} = [n]_q,
\quad \qbi{n}{p, 2}{q} = \frac{[n]_q}{[2]_q}([n-1]_q+[p]_q[n-p]_q).\]

Applying \eqref{q-whip-2} to the definition \eqref{eq:def2}, we get
\[\qbi{n}{p,k}{q} = \sum_{r \ge 0} (-1)^i
\qbi{n-i}{k}{q}\frac{[n]}{[n-i]}\qbi{n-i}{i}{q} \qbi{n-2i}{p-i}{q}
q^{i(i-1)/2+ki}.\label{eq:expression1}
\]

The following two generating functions are obtained in the same
way as for ${n\choose p,k}_q$.
\begin{thm} There holds
\begin{eqnarray}
& &\sum_{0 \leq k,p \leq n}{n \brack p,k}_qx^p q^{p \choose 2} y^k
q^{k \choose 2} \nonumber\\[5pt]
&=& (-y;q)_n(-x;q)_n\sum_{i \geq 0} \frac{1-q^n}{1-q^{n-i}}{n-i
\brack
i}_q\frac{(-x)^iq^{i^2-i}}{(-y;q)_i(-x;q)_{i}(-xq^{n-i};q)_i}.\label{gf2}
\end{eqnarray}
\end{thm}

\begin{thm} For $0\leq p\leq n$, there holds
\begin{eqnarray}
\sum_{k \geq 1}{n \brack p,k}_qy^k q^{{k \choose 2}} =
y(-yq^{n-p};q)_p[n]_q\,{}_2\phi_1 \!
\left[\left.\begin{matrix} q^{p+1},q^{p-n+1}\\
q^2 \end{matrix}\right|q;-yq^{n-p}\right].\label{eq:gf1}
\end{eqnarray}
\end{thm}

As \eqref{conti2} we get the recurrence relation
\begin{prop} For $k\neq 0$ and $0\leq p\leq n$, there holds
\begin{eqnarray}
(1-q^{n-p+1}){n \brack p-1, k}_q -(1-q^p){n \brack p, k}_q =
\frac{[n]_q}{[n-1]_q}(1-q^{n-2p+1})q^{k+p-1}{n-1 \brack p-1, k}_q.
\end{eqnarray}
\end{prop}

\begin{cor} We have
\[ {n \brack p,k}_q=\frac{[n]_q}{[p]_q}\,\sum_{i=0}^{k-1}{n-p\brack
k-1-i}_q{n-p+i\brack i}_q{p\brack i+1}_qq^{(i+1+p-n)(i+1-k)}.
\]
\end{cor}
As in the proof of Theorem 3, we can prove that ${n \brack p,k}_q$
is also a q-positive integer by using the above Corollary.

\section{Tables of the generalized $q$-binomial coefficients}

 \begin{itemize}
 \item ${n \brack 0, k}_q = {n \choose 0,k}_q = {n \brack k}_q$ for $k \geq 0$,
 \item ${n \brack p, 0}_q = {n \choose p,0}_q = 0$, for $p \geq  0$.
 \item $n=1, {n \brack 1,1}_q={n \choose 1,1}_q = 1$.
 \end{itemize}

\begin{center}
{\bf Tables of ${n \choose p,k}_q$}
\end{center}
\[
\begin{array}{ll}
n=2 & \footnotesize{
\begin{tabular}{|c|c||c|c|}
\hline  $p \backslash k $ & 1 & 2\\ \hline \hline 1 & $(1+q)q$ &
1+q
\\\hline
 2 & $1+q$ & 1 \\\hline
\end{tabular}}\\ \\
n=3, &{\footnotesize{
\begin{tabular}{|c||c|c|c|}
 \hline $p \backslash k$  & 1 & 2 & 3 \\\hline \hline 1  & $q^2[3]_q $ &
$q(1+q)(1+q+q^2)$ & $[3]_q$
\\ \hline
 2 & $ q^2[3]_q$ & $q(1+q)(1+q+q^2)$ &  $[3]_q$\\ \hline
 3 & $[3]_q$ & $[3]_q$ & 1 \\ \hline
\end{tabular}}}\\ \\
$n=4$, & \footnotesize{\begin{tabular}{|c||c|c|c|c|}
\hline $p \backslash k$  & 1 & 2 & 3 & 4 \\ \hline \hline 1  &
$q^3[4]_q$ & $q^2[3]_q[4]_q$ & $q[2]_q[3]_q(1+q^2)$ & $[4]_q$
\\\hline
 2 & $q^4[4]_q$ & $(1+q^2)(1+3q+2q^2+q^3)$ &  $[2]_q(1+q^2)(2+q+q^2)$ & $[3]_q(1+q^2)$\\\hline
 3 & $q^3[4]_q$ & $q^2[3]_q[4]_q$ & $q[2]_q[3]_q(1+q^2)$ & $[4]_q$\\\hline
 4 & $[4]_q$ & $[2]_q[3]_q$ & $[4]_q$ & 1 \\ \hline
\end{tabular}}\\ \\
$n=5$, & \footnotesize{
\begin{tabular}{|c||c|c|c|c|c|}
 \hline $p \setminus k$ & 1 & 2 & 3 & 4 & 5\\  \hline \hline
1 & $[5]_q$ & $[2]_q[5]_q(1+q^2)$ & $[3]_q[5]_q(1+q^2)$ & $[2]_q[5]_q(1+q^2)$ & $[5]_q$  \\  \hline
2 & $[5]_q$ & $[5]_q(1+q+2q^2+q^3)$ & $[3]_q[5]_q(1+2q^2)$ & $[5]_q(1+q+2q^2+2q^3+q^4)$ & $[5]_q(1+q^2)$ \\  \hline
3 & $[5]_q$ & $[5]_q(1+q+2q^2+q^3)$ & $[3]_q[5]_q(1+2q^2)$ &  $[5]_q(1+q+2q^2+2q^3+q^4)$ & $[5]_q(1+q^2)$ \\  \hline
4 & $[5]_q$ & $[2]_q[5]_q(1+q^2)$ & $[3]_q[5]_q(1+q^2)$ & $[2]_q[5]_q(1+q^2)$ & $[5]_q$\\  \hline
5 & $[5]_q$ & $[5]_q(1+q^2)$ & $[5]_q(1+q^2)$ & $[5]_q$ & 1 \\ \hline
\end{tabular}}
\end{array}
\]

\newpage
\begin{center}
{\bf Tables of ${n \brack p, k}_q$}
\end{center}

\[
\begin{array}{ll}
n=2, & \footnotesize{\begin{tabular}{|c||c|c|}
 \hline $p \backslash k $ & 1 & 2\\ \hline \hline
1 & $[2]_q$ & 1+q \\\hline
 2 & $[2]_q$ & 1 \\\hline
\end{tabular}}\\ \\

n=3, & \footnotesize{\begin{tabular}{|c||c|c|c|}
\hline $p \backslash k$  & 1 & 2 & 3 \\  \hline \hline 1  &
$[3]_q$ & $[2]_q[3]_q$ & $[3]_q$
\\\hline
 2 & $[3]_q$ & $[2]_q[3]_q $ &  $[3]_q$\\\hline
 3 & $[3]_q$ & $[3]_q $ & 1\\\hline
\end{tabular}}\\ \\

n=4, &\footnotesize{\begin{tabular}{|c||c|c|c|c|}
\hline $p \backslash k$  & 1 & 2 & 3 & 4 \\\hline \hline 1  &
$[4]_q$ & $[3]_q[4]_q$ & $[2]_q[3]_q(1+q^2)$ &$[4]_q$
\\\hline
 2 & $[4]_q$ & $(1+q^2)(1+2q+3q^2+q^3)$ &  $[2]_q(1+q^2)(1+q+2q^2)$ & $(1+q^2)[3]_q$\\\hline
 3 & $[4]_q$ & $[3]_q[4]_q$ & $[2]_q[3]_q(1+q^2)$ & $[4]_q$\\\hline
 4 & $[4]_q$ & $[3]_q(1+q^2)$ & $[4]_q$ & 1 \\\hline
\end{tabular}}\\ \\

n=5, & \footnotesize{
\begin{tabular}{|c||c|c|c|c|c|}
\hline $p \setminus k$ & 1 & 2 & 3 & 4 & 5\\ \hline \hline
1 & $q^4[5]_q$ & $[2]_q[5]_q(1+q^2)q^3$ & $[3]_q[5]_q(1+q^2)q^2$ & $[2]_q[5]_q(1+q^2)q$  & $[5]_q$ \\\hline
2 & $q^6[5]_q$ & $[5]_q(1+2q+q^2+q^3)q^5$ & $[3]_q[5]_q(2+q^2)q^4$ & $[5]_q(1+2q+2q^2+q^3+q^4)q^2$ & $[5]_q(1+q^2)$\\\hline
3 & $q^6[5]_q$ & $[5]_q(1+2q+q^2+q^3)q^5$ & $[3]_q[5]_q(2+q^2)q^4$ & $[5]_q(1+2q+2q^2+q^3+q^4)q^2$ & $[5]_q(1+q^2)$\\\hline
4 & $q^4[5]_q$ & $[2]_q[5]_q(1+q^2)q^3$ & $[3]_q[5]_q(1+q^2)q^2$ & $[2]_q[5]_q(1+q^2)q$ & $[5]_q$ \\\hline
5 & $[5]_q$ & $[5]_q(1+q^2)$ & $[5]_q(1+q^2)$ & $[5]_q$ & 1 \\ \hline
\end{tabular}}
\end{array}
\]

From the above tables, it seems that the sequences of the
coefficients in ${n \choose p, k}$ and ${n \brack p, k}$ are
unimodal for each fixed pair $(n, p)$ or $(n,k)$.

\section*{Acknowledgements} This work was done under the auspices of
the National Science Foundation of China.

\newpage


\begin{thebibliography}{9}
\bibitem{GR}
   G. Gasper \& M. Rahman,
\emph{Basic hypergeometric series},
   Encyclopedia of Math. and its Applications, \textbf{35} (1990).
 \bibitem{JLZ}
   F. Jouhet, B. Lass and  J. Zeng,
   \emph{Sur une g\'en\'eralisation des coefficients binomiaux},
   preprint, arXiv:math.CO/0303025 v1, 3 Mars 2003.
   \bibitem{La}
   M. Lassalle, \emph{A new family of positive integers},
   Ann. Comb. 6(2002), no. 3-4, 399-405.
   \bibitem{La2}
   M. Lassalle, \emph{Jack polynomials and some identities for
   partitions}, to appear in Trans. of Amer. Math. Soc., 2004.
\end{thebibliography}
\end{document}